\documentclass[11pt]{article}

\usepackage[margin=1in]{geometry}
\usepackage{amsmath,amssymb,amsthm,mathtools}
\usepackage[mathscr]{euscript}
\usepackage[T1]{fontenc}
\usepackage[utf8]{inputenc}
\usepackage{lmodern}
\usepackage{enumitem}
\usepackage{microtype}
\usepackage[nocompress]{cite}
\usepackage[colorlinks=true,linkcolor=blue,citecolor=blue,urlcolor=blue]{hyperref}
\hypersetup{
  pdftitle={Recursive Parameter Recovery for the Modified Leray-Alpha Model via Continuous Data Assimilation},
  pdfauthor={},
  pdfsubject={Continuous data assimilation and parameter recovery},
  pdfkeywords={modified Leray-alpha model, continuous data assimilation, parameter recovery, determining modes, nudging}
}

\theoremstyle{plain}
\newtheorem{theorem}{Theorem}
\newtheorem{proposition}{Proposition}
\newtheorem{lemma}{Lemma}
\theoremstyle{definition}
\newtheorem{assumption}{Assumption}
\newtheorem{definition}{Definition}
\newtheorem{algorithm}{Algorithm}
\theoremstyle{remark}
\newtheorem{remark}{Remark}
\theoremstyle{plain}

\newcommand{\norm}[1]{\left\|#1\right\|}

\title{Parameter Recovery for the Modified Leray-\texorpdfstring{$\alpha$}{alpha} Model via Continuous Data Assimilation}

\author{D\'ebora A.F. Albanez\textsuperscript{a}\thanks{Corresponding author: D\'ebora A.F. Albanez. E-mail addresses: deboraalbanez@utfpr.edu.br (D.A.F. Albanez), mjb26b@acu.edu \& m.benvenutti@ufsc.br (M.J. Benvenutti), jtian@towson.edu (J. Tian).}\quad Maicon J. Benvenutti\textsuperscript{b, c}\quad Jing Tian\textsuperscript{d}\\[2mm]
\begin{minipage}{0.92\textwidth}\centering\small
\textsuperscript{a}Departamento Acad\^emico de Matem\'atica, Universidade Tecnol\'ogica Federal do Paran\'a,\\
Corn\'elio Proc\'opio PR, 86300-000, Brazil\\
\textsuperscript{b} Abilene Christian University, Abilene, TX 79699, USA\\
\textsuperscript{c} 
Universidade Federal de Santa Catarina, Blumenau SC 89065-200, Brazil\\
\textsuperscript{d}Department of Mathematics, Towson University, Towson, MD 21252, USA
\end{minipage}}
\date{}

\begin{document}
\maketitle

\begin{abstract}
In this study, we extend the recursive length-scale recovery method for the simplified Bardina model to the three-dimensional modified Leray-$\alpha$ model. The main difficulty is that the unknown parameter also occurs in the nonlinear convection term, so estimating the error in the parameter update requires control of second spatial derivatives of the velocity error. We combine continuous data assimilation with a recursive update of the squared length scale, using observed low Fourier modes and computed high Fourier modes. Under uniform regularity bounds on the reference solution and suitable conditions on the observations and the recovery procedure, we prove geometric convergence of both the parameter and velocity errors at the update times. When the update-window lengths are uniformly bounded, this convergence is exponential in time at the update times.
\end{abstract}



\section{Introduction}\label{sec:introduction}

The Navier--Stokes equations describe the motion of viscous incompressible fluids. In turbulent regimes, however, the range of dynamically relevant scales makes direct numerical simulation expensive. Regularized fluid models provide a way to describe the large-scale motion without resolving every small scale. Beginning with Leray's regularization \cite{Leray1934}, this approach has led to several families of turbulence models, including the Navier--Stokes-$\alpha$, Leray-$\alpha$, modified Leray-$\alpha$, Clark-$\alpha$, and simplified Bardina models \cite{FoiasHolmTiti2001,CheskidovHolmOlsonTiti2005,IlyinLunasinTiti2006,CaoHolmTiti2005,Bardina1980,CaoLunasinTiti2006,LaytonLewandowski2006}. The parameter $\alpha$ specifies a filtering length, and its choice influences the relation between the resolved velocity and the underlying regularized flow. Recovering this parameter from observations is therefore a natural problem in model calibration.

In this paper, we consider the three-dimensional modified Leray-$\alpha$ model introduced in \cite{IlyinLunasinTiti2006}. On the periodic box $\Omega=[0,L_\Omega]^3$, its equations are
\begin{equation}\label{eq:MLalpha-physical}
\begin{cases}
v_t-\nu\Delta v+(v\cdot\nabla)u+\nabla p=f,\\
\nabla\cdot u=\nabla\cdot v=0,\\
v=u-\alpha^2\Delta u.
\end{cases}
\end{equation}
Here $u$ is the filtered velocity, $v$ is the associated unfiltered velocity, $p$ is the pressure, $\nu>0$ is the viscosity, and $f$ is the external force. 
We impose periodic boundary conditions and assume
\[
u(\cdot,0)=u_0,\qquad
\nabla\cdot f=0,\qquad
\int_\Omega u_0(x)\,dx
=\int_\Omega f(x,t)\,dx=0,
\quad t\ge0.
\]
The force and viscosity are known, whereas the initial state and the parameter $\alpha$ are unknown. Only a finite number of low Fourier modes of $u$ are observed. We assume that known constants $\alpha_0$ and $\alpha_1$ satisfy
\[
0<\alpha_0\le\alpha\le\alpha_1.
\]

The theory of determining modes and nodes shows that sufficiently resolved partial observations determine the asymptotic behavior of the full solution under model-specific hypotheses \cite{FoiasProdi1967,FoiasTemam1984,JonesTiti1993}. Continuous data assimilation incorporates these observations into an evolution equation through a feedback term. We use the framework introduced by Azouani, Olson, and Titi \cite{AOT2014}, whose applications to turbulence models include \cite{AlbanezLopesTiti2016,AlbanezBenvenutti2018}. When some parameters are also unknown, the feedback system can be coupled with a parameter-recovery procedure. Such methods have been studied for viscosity recovery and sensitivity analysis \cite{CarlsonHudsonLarios2020,CarlsonLarios2021,Martinez2022,BiswasHudson2023}, identification of forcing \cite{FarhatLariosMartinezWhitehead2024}, and simultaneous recovery of multiple parameters \cite{CarlsonEtAl2022,PachevWhiteheadMcQuarrie2022}. Related developments include on-the-fly model discovery and relaxation-based estimation \cite{NeweyWhiteheadCarlson2025,MartinezMurriWhitehead2025}, viscosity recovery for the steady Navier--Stokes equations \cite{RebholzReyesWhitehead2026}, recovery for the primitive equations with a magnetic field \cite{ZhaoEtAl2026}, and data-driven operator inference \cite{WangYan2026}.

For $\alpha$-regularization models, the effect of a mismatch in the filtering length has been investigated in \cite{ParameterAnalysis2025,MLAlphaParameter2025}. A recursive recovery procedure for the simplified Bardina model was subsequently developed in \cite{BardinaParameter2026}. We adapt that procedure to recover $\alpha$ and approximate the reference solution of the modified Leray-$\alpha$ equation. Unlike the simplified Bardina model, whose convection term is $(u\cdot\nabla)u$, \eqref{eq:MLalpha-physical} contains $(v\cdot\nabla)u$. When the system is written in terms of $u$, the unknown parameter $\alpha$ therefore occurs not only in the time-derivative and viscous terms, but also in the convection term. To construct the approximate solution, let $P_N$ denote the orthogonal projection onto the spatial Fourier modes indexed by $k\in\mathbb Z^3$ with $0<|k|\le N$. For a fixed approximate parameter $\beta\in[\alpha_0,\alpha_1]$, we use the same feedback system as in \cite{MLAlphaParameter2025} with Fourier-mode observations:
\begin{equation}\label{eq:assim-physical}
\begin{cases}
z_t-\nu\Delta z+(z\cdot\nabla)w+\nabla q
=f-\eta(I-\beta^2\Delta)P_N(w-u),\\
\nabla\cdot w=\nabla\cdot z=0,\\
z=w-\beta^2\Delta w.
\end{cases}
\end{equation}
The fields $w$ and $z$ satisfy the same periodic and zero-mean conditions as $u$ and $v$; $q$ is the corresponding pressure, $I$ is the identity, and $\eta>0$ is the nudging strength. The initial value of $w$ is arbitrary in the energy space defined below. 

The parameter update depends on the unknown velocity, whereas only its low Fourier modes are observed. Following \cite{BardinaParameter2026}, we replace the unknown velocity with a field formed from the observed low Fourier modes and the high Fourier modes of the approximate solution. For the modified Leray-$\alpha$ model, this replacement introduces an error involving second spatial derivatives, which cannot be controlled by a basic velocity error estimate alone. We combine a synchronization estimate for the velocity error with a time-integrated higher-order dissipation estimate to control the second spatial derivatives appearing in the parameter update. These estimates are then used to construct a recursive parameter update. Before each update window, we keep the parameter fixed for a relaxation period to allow the assimilated solution to approach the reference solution. Under suitable regularity assumptions on the reference solution and sufficient observational signal relative to the approximation error, we prove that both the estimated parameter and the approximate velocity converge geometrically to their corresponding true values at the update times. Furthermore, if the update intervals are uniformly bounded in length, this convergence is exponential in time at the update times. 

The rest of the paper is organized as follows. Section~\ref{sec:functional} presents the functional setting. Section~\ref{sec:algorithm} introduces the parameter recovery procedure. Section~\ref{sec:reference} establishes the estimates for the reference and data-assimilated systems. Section~\ref{sec:convergence} presents and proves the main theorem. Section~\ref{sec:conclusion} concludes the paper.

\section{Functional analytic framework}\label{sec:functional}

For $\Omega=[0,L_\Omega]^3$, we denote by $L^p(\Omega)$ the usual Lebesgue spaces of three-dimensional vector fields. For each $s\in\mathbb R$, we define the Hilbert space
\[
\dot V_s=\left\{
\begin{aligned}
\phi&=\sum_{k\in\mathbb Z^3\setminus\{0\}}
\hat\phi_k e^{2\pi i k\cdot x/L_\Omega}:\quad
\hat\phi_k\in\mathbb C^3,\\
&\hat\phi_{-k}=\overline{\hat\phi_k},\quad k\cdot\hat\phi_k=0,\quad
\sum_{k\in\mathbb Z^3\setminus\{0\}}
\left(\frac{2\pi|k|}{L_\Omega}\right)^{2s}
|\hat\phi_k|^2<\infty
\end{aligned}
\right\},
\]
where $\hat\phi_k$ are the Fourier coefficients of $\phi$, and the series is understood in the sense of distributions.
We endow this space with the norm
\[
\norm{\phi}_{\dot V_s}^2
=|\Omega|\sum_{k\ne0}
\left(\frac{2\pi|k|}{L_\Omega}\right)^{2s}|\hat\phi_k|^2.
\]
So $\norm{\phi}_{\dot V_0}=\norm{\phi}_{L^2}$.
We write
\[
H=\dot V_0,\quad V=\dot V_1,\quad V'=\dot V_{-1},\quad D'=\dot V_{-2},\quad \dot V_{-s}=(\dot V_s)'\quad(s\ge0).
\]

We denote by $\mathscr P:L^2(\Omega)\to H$ the
Helmholtz--Leray orthogonal projection, given by
\[
\mathscr P\phi
=\sum_{k\in\mathbb Z^3\setminus\{0\}}
\left(\hat\phi_k-
\frac{k(k\cdot\hat\phi_k)}{|k|^2}\right)
e^{2\pi i k\cdot x/L_\Omega}.
\]
The same formula defines its extension to periodic distributions.
For each $s\in\mathbb R$, the Stokes operator
$A:\dot V_s\to\dot V_{s-2}$ is given by
\[
A\phi
=\sum_{k\in\mathbb Z^3\setminus\{0\}}
\frac{4\pi^2|k|^2}{L_\Omega^2}\,
\hat\phi_k e^{2\pi i k\cdot x/L_\Omega}.
\]
In particular,
\[
A\phi=-\mathscr P\Delta\phi=-\Delta\phi
\quad(\phi\in D(A)=\dot V_2),
\qquad \lambda_1=\frac{4\pi^2}{L_\Omega^2}.
\]

For the inner product and norms, we use the notation
\[
(\phi,\psi)=(\phi,\psi)_{L^2},\qquad
\norm{\phi}=\norm{\phi}_{L^2},\qquad
\norm{\phi}_V=\norm{\nabla\phi}.
\]
Thus round parentheses denote the $L^2$ inner product, whereas $\langle\cdot,\cdot\rangle$
denotes the appropriate Sobolev duality pairing, which extends this inner product.
The zero-mean condition gives the Poincar\'e inequalities
\begin{equation}\label{eq:poincare}
\norm{u}^2\le\lambda_1^{-1}\norm{\nabla u}^2\quad(u\in V),\qquad
\norm{\nabla u}^2\le\lambda_1^{-1}\norm{Au}^2\quad(u\in D(A)).
\end{equation}
We recall the following particular cases of the Gagliardo--Nirenberg inequalities \cite{FoiasManleyRosaTemam2001,Temam1995}:
\begin{equation*}
\norm{g}_{L^3}\le C\norm{g}^{1/2}\norm{\nabla g}^{1/2},\qquad
\norm{g}_{L^4}\le C\norm{g}^{1/4}\norm{\nabla g}^{3/4},\qquad
\norm{g}_{L^6}\le C\norm{\nabla g}.
\end{equation*}
These inequalities hold for every $g\in V$.
We shall also use the Sobolev embeddings $\dot V_2\hookrightarrow L^\infty$ and
$\dot V_3\hookrightarrow W^{1,\infty}$.

For each $\vartheta>0$, we denote the Helmholtz operator by $T_\vartheta=I+\vartheta^2A$.
Its inverse is the Fourier multiplier operator with multiplier
$(1+\vartheta^2\lambda_1|k|^2)^{-1}$, and the following estimates hold:
\begin{equation}\label{eq:filter-inverse}
\norm{T_\vartheta^{-1}\phi}\le\norm{\phi}\quad(\phi\in H),\qquad
\norm{T_\vartheta^{-1}\mathcal F}\le
\vartheta^{-2}\norm{\mathcal F}_{D'}
\quad(\mathcal F\in D').
\end{equation}
The observations are represented by the low-mode projection introduced in the Introduction. For $N\in\mathbb N$, this projection and its complement are given by
\[
P_N\phi=\sum_{0<|k|\le N}\hat\phi_ke^{2\pi i k\cdot x/L_\Omega},
\qquad Q_N=I-P_N.
\]
By orthogonality in each space $\dot V_s$, these projections satisfy the following relations:
\begin{equation}\label{eq:projection-properties}
\norm{Q_N\phi}^2\le\frac{1}{\lambda_1N^2}\norm{\nabla \phi}^2,
\qquad \norm{P_N\phi}\le\norm{\phi},\qquad
\norm{\nabla P_N\phi}\le\norm{\nabla \phi},\qquad P_NA=AP_N.
\end{equation}

We define the bilinear operator $B:V\times V\to V'$ by
\[
B(u,v)=\mathscr P[(u\cdot\nabla)v].
\]
It extends continuously to $B:H\times V\to D'$, with 
\begin{equation}\label{eq:negative-bilinear}
\norm{B(u,v)}_{D'}\le C\norm{u}\norm{\nabla v}\qquad(u\in H,\ v\in V).
\end{equation}
For $u,v,w\in V$, and also for $u\in H$ and $v,w\in D(A)$, we have
\begin{equation}\label{eq:skew-symmetry}
\langle B(u,v),v\rangle=0,\qquad
\langle B(u,v),w\rangle=-\langle B(u,w),v\rangle.
\end{equation}
Throughout the paper, $C$ denotes a positive constant that may
change from one occurrence to another, and may depend on the domain.

\section{Parameter recovery algorithm}\label{sec:algorithm}

Throughout the analysis, we assume that $f\in L^\infty(0,\infty;H)$ and write
$F=\norm{f}_{L^\infty(0,\infty;H)}$.
Applying $\mathscr P$ to \eqref{eq:MLalpha-physical} gives the functional formulation
\begin{equation}\label{eq:MLalpha-v}
\begin{cases}
\displaystyle \frac{dv}{dt}+\nu Av+B(v,u)=f,\\[1mm]
v=T_\alpha u=(I+\alpha^2A)u,\\
u(0)=u_0.
\end{cases}
\end{equation}

Similarly, for fixed $\beta\in[\alpha_0,\alpha_1]$, the assimilation system \eqref{eq:assim-physical} becomes
\begin{equation}\label{eq:assim-functional}
\begin{cases}
\displaystyle \frac{dz}{dt}+\nu Az+B(z,w)
=f-\eta T_\beta P_N(w-u),\\[1mm]
z=T_\beta w=(I+\beta^2A)w.
\end{cases}
\end{equation}
Here the feedback field $P_N(w-u)$ is divergence free. 

For sufficiently regular $\phi$, we introduce two operators that separate the terms with and without the squared filtering length:
\begin{equation*}
\begin{aligned}
R(\phi)&=\phi_t+\nu A\phi+B(\phi,\phi),\\
S(\phi)&=A\phi_t+\nu A^2\phi+B(A\phi,\phi).
\end{aligned}
\end{equation*}
With this notation, systems \eqref{eq:MLalpha-v} and \eqref{eq:assim-functional} can be written in the form
\begin{equation}\label{eq:functional-RS}
R(u)+\alpha^2S(u)=f,\qquad
R(w)+\beta^2S(w)=f-\eta T_\beta P_N(w-u).
\end{equation}
Estimate \eqref{eq:negative-bilinear} yields 
\[
\phi\in L^\infty(V)\cap L^2(D(A)),\quad \phi_t\in L^2(H)
\quad\Longrightarrow\quad
A\phi_t,A^2\phi,B(A\phi,\phi)\in L^2(D').
\]
Fixing $t_1\ge0$, a regular solution satisfies
\[
u\in C(\mathcal J;V)\cap L^2(\mathcal J;D(A)),\qquad
u_t\in L^2(\mathcal J;H),
\]
on a bounded interval $\mathcal J=[r,r+\ell]$, where
$r\ge t_1$ and $0<\ell<\infty$.

For $\vartheta>0$ and $\phi\in V$, we define the corresponding energy by
\begin{equation*}
E_\vartheta(\phi)=\norm{\phi}^2+\vartheta^2\norm{\nabla \phi}^2.
\end{equation*}
For $u_0\in V$, the reference system admits a unique global regular solution:
see \cite{IlyinLunasinTiti2006} for autonomous forcing and
\cite{MLAlphaParameter2025} for the time-dependent forcing considered here. Taking the duality pairing of \eqref{eq:MLalpha-v}
with $u$ and using \eqref{eq:skew-symmetry} and Young's inequality together with
\eqref{eq:poincare}, we obtain
\[
\frac{d}{dt}E_\alpha(u)+\nu\lambda_1E_\alpha(u)
\le\frac{F^2}{\nu\lambda_1}.
\]
Applying Gronwall's inequality, we conclude that, for $t\ge t_1$,
\begin{equation}\label{eq:reference-energy}
E_\alpha(u(t))\le e^{-\nu\lambda_1(t-t_1)}E_\alpha(u(t_1))
+\frac{F^2}{\nu^2\lambda_1^2}
\bigl(1-e^{-\nu\lambda_1(t-t_1)}\bigr).
\end{equation}
Assume that a known, parameter-independent bound
$E_{\alpha_1}(u(t_1))^{1/2}\le R_{\rm ref}$ is available.
Then \eqref{eq:reference-energy} gives
\begin{equation}\label{eq:Mref}
E_\alpha(u(t))^{1/2}\le
M_{\rm ref}:=\max\left\{R_{\rm ref},\frac{F}{\nu\lambda_1}\right\},
\qquad t\ge t_1.
\end{equation}
For the parameter recovery estimates, we impose the following additional regularity assumptions on the reference solution.

\begin{assumption}\label{ass:reg}
There exist constants $M_3,M_4>0$, depending only on the
reference solution, such that
\textup{\textbf{(i)}}
\begin{equation}\label{eq:reg1}
\sup_{t\ge t_1}\norm{u(t)}_{\dot V_3}\le M_3.
\end{equation}
\textup{\textbf{(ii)}}
\begin{equation*}
\sup_{r\ge t_1}\int_r^{r+1}\left(\norm{u(s)}_{\dot V_4}^2+\norm{u_t(s)}_{\dot V_2}^2\right)ds\le M_4^2.
\end{equation*}
\end{assumption}
For every $T>t_1$, condition~\textup{(ii)} and the embedding
$L^2(t_1,T;\dot V_4)\cap H^1(t_1,T;\dot V_2)
\hookrightarrow C([t_1,T];\dot V_3)$ give a continuous representative
of $u$, which we use throughout. Thus the supremum in
\eqref{eq:reg1} agrees with the essential supremum.

As in \cite{BardinaParameter2026}, we keep $N$ and $\eta$ fixed throughout the recovery iteration. Let $t_1<t_2<\cdots$ denote the successive parameter-update times,
with $\beta_n$ held fixed on $[t_n,t_{n+1}]$.
On each interval $[t_n,t_{n+1}]$, we use the approximate parameter $\beta_n\in[\alpha_0,\alpha_1]$ and solve
\begin{equation}\label{eq:assim-system}
\begin{cases}
(z_n)_t+\nu Az_n+B(z_n,w_n)=f-\eta T_{\beta_n}P_N(w_n-u),\quad z_n=T_{\beta_n}w_n,\\[1mm]
w_n(t_n)=w_{n-1}(t_n)\ (n\ge2),\qquad w_1(t_1)=w_0.
\end{cases}
\end{equation}
The initial state $w_0\in V$ is arbitrary. 

Let $g_n=w_n-u$, so that $P_Ng_n=P_Nw_n-P_Nu$.
Subtracting the two equations in \eqref{eq:functional-RS}, we obtain
\begin{equation}\label{eq:key-identity}
(\alpha^2-\beta_n^2)S(u)=R(w_n)-R(u)+\beta_n^2\big(S(w_n)-S(u)\big)+\eta(I+\beta_n^2A)P_Ng_n.
\end{equation}

Following \cite{BardinaParameter2026}, we approximate the reference solution by combining the observed low modes of $u$ with the computed high modes of $w_n$:
\begin{equation}\label{eq:usharp}
u_N^{*,n}=w_n-P_Ng_n=Q_Nw_n+P_Nu,
\qquad u-u_N^{*,n}=-Q_Ng_n.
\end{equation}
Define
\begin{equation}\label{eq:Dn}
D_N^{(n)}(t):=A(P_Nu)_t+\nu A^2P_Nu+P_NB(Au_N^{*,n},u_N^{*,n}).
\end{equation}
Using the same replacement in the projected error equation, we define
\begin{equation}\label{eq:Yn}
\begin{aligned}
Y_N^{(n)}:={}&T_{\beta_n}\bigl(P_N(g_n)_t+\nu AP_Ng_n+\eta P_Ng_n\bigr)\\
&+P_N\bigl\{B(T_{\beta_n}w_n,w_n)-B(T_{\beta_n}u_N^{*,n},u_N^{*,n})\bigr\}.
\end{aligned}
\end{equation}
Since $P_NA=AP_N$, \eqref{eq:Dn} and \eqref{eq:Yn} can be computed from
$w_n$, $P_Nu$, and $(P_Nu)_t$. Moreover, 
\[
D_N^{(n)},Y_N^{(n)}\in L^2(t_n,t_{n+1};P_NH).
\]
We use exact observations in continuous time.

\begin{definition}\label{def:windows}
For each $n\ge1$, we choose an intermediate time $\widehat t_n$ and denote the corresponding update interval and its length by
\begin{equation*}
t_n<\widehat t_n<t_{n+1},\qquad I_n=[\widehat t_n,t_{n+1}],\qquad \tau_n=|I_n|.
\end{equation*}
The interval $[t_n,\widehat t_n]$ is used for relaxation, and $I_n$ is used to compute the parameter update. We require the update signal to satisfy
\begin{equation}\label{eq:delta}
\delta_n:=\int_{\widehat t_n}^{t_{n+1}}\norm{D_N^{(n)}(s)}^2ds>0.
\end{equation}
Fix $\zeta>0$. An update interval $I_n$ is accepted if
$\delta_n/\tau_n\ge\zeta$.
\end{definition}

On each update interval, we first form a candidate value for the squared parameter:
\begin{equation}\label{eq:update}
\gamma_{n+1}:=\beta_n^2+\frac{1}{\delta_n}\int_{\widehat t_n}^{t_{n+1}}\big(Y_N^{(n)}(s),D_N^{(n)}(s)\big)ds,
\end{equation}
then project this value onto $[\alpha_0^2,\alpha_1^2]$, obtaining
\begin{equation}\label{eq:clipping}
\beta_{n+1}
=\sqrt{\min\!\left\{\alpha_1^2,\,
\max\!\left\{\alpha_0^2,\gamma_{n+1}\right\}\right\}}.
\end{equation}

We propose the following parameter updating algorithm.

\begin{algorithm}\label{alg:recovery}
\mbox{}\par\smallskip
\noindent\textbf{Initialization.}
Given the observations $P_Nu(t)$ on $[t_1,T]$, where $T>t_1$,
choose $w_0\in V$ and $\beta_1\in[\alpha_0,\alpha_1]$.
Let $L>0$ be the fixed relaxation time before each update window.
Keep $N\in\mathbb N$, $\eta>0$, and the threshold $\zeta>0$ fixed, and choose a finite list of trial lengths satisfying
\[
0<\tau_0=\ell_1<\cdots<\ell_J\le\tau_{\max},\qquad J\ge1.
\]
Require the reference estimates from $t_1$ onward; set $w_1(t_1)=w_0$, $n=1$.

\smallskip
\noindent\textbf{Iteration.}
Given $t_n,\beta_n,w_n(t_n)$:
\begin{enumerate}[label=\arabic*.,leftmargin=2em,itemsep=4pt,topsep=4pt]
\item If $t_n+L+\tau_0>T$, stop. Otherwise set $\widehat t_n=t_n+L$.
\item Try the candidate window lengths $\ell_1,\ldots,\ell_J$ in increasing order, restricted to $\widehat t_n+\ell_j\le T$. Set
\[
\tau_n=\ell_j,\qquad t_{n+1}=\widehat t_n+\tau_n.
\]
On $[t_n,t_{n+1}]$, solve \eqref{eq:assim-system} with the same $\beta_n,w_n(t_n)$ for every trial.
Compute \eqref{eq:usharp}--\eqref{eq:delta}; accept the first trial with
$\delta_n/\tau_n\ge\zeta$.
\item If all eligible trials fail, stop without changing $\beta_n$ or $w_n(t_n)$.
\item For an accepted window, apply \eqref{eq:update}--\eqref{eq:clipping}, set
$w_{n+1}(t_{n+1})=w_n(t_{n+1})$, increment $n$, and repeat.
\end{enumerate}
\noindent\textbf{Output.}
Return $\beta_n,w_n(t_n)$, or $\beta_1,w_0$ if no update was accepted.
\end{algorithm}

\section{Estimates for the reference and assimilation systems}\label{sec:reference}

For each $\kappa>0$, we introduce the quantities
\begin{equation}\label{eq:reference-constants}
\begin{aligned}
K_S(\kappa)&=C\left(\frac{(1+\nu^2)M_4^2}{1-e^{-\kappa}}+\frac{M_3^4}{\kappa}\right),\\
M_S^2&=C\left((1+\nu^2)M_4^2(1+\tau_0^{-1})+M_3^4\right),
\end{aligned}
\end{equation}
where $\tau_0$ is the shortest trial length and $C$ is a sufficiently large
constant depending only on the domain.

\begin{lemma}
Under Assumption~\ref{ass:reg}, for $\mathcal J=[r,r+\ell]\subset[t_1,\infty)$,
with $0<\ell<\infty$, we have $S(u)\in L^2(\mathcal J;H)$ and
\begin{equation}\label{eq:S-control}
\int_{\mathcal J}\norm{S(u(t))}^2\,dt
\le C\bigl((1+\nu^2)M_4^2(\ell+1)+M_3^4\ell\bigr)
\le M_S^2\ell\qquad(\ell\ge\tau_0>0).
\end{equation}
 Moreover, for every $\kappa>0$,
\begin{equation}\label{eq:S-weighted}
\int_r^t e^{-\kappa(t-s)}\norm{S(u(s))}^2\,ds
\le K_S(\kappa),\qquad t\ge r\ge t_1.
\end{equation}
\end{lemma}

\begin{proof}
Applying Assumption~\ref{ass:reg} on successive unit intervals
starting at $r$, we obtain
\begin{equation}\label{eq:unit-bound}
\int_{\mathcal J}\left(\norm{u_t}^2_{\dot V_2}+\norm{u}^2_{\dot V_4}\right)ds\le M_4^2(\ell+1).
\end{equation}
For the nonlinear term, H\"older's inequality and the Sobolev
embeddings give, for a.e. $s\ge t_1$,
\[
\norm{B(Au(s),u(s))}
\le C\norm{Au(s)}_{L^3}\norm{\nabla u(s)}_{L^6}
\le C\norm{u(s)}_{\dot V_3}^2
\le CM_3^2.
\]
Since $S(u)=Au_t+\nu A^2u+B(Au,u)$, we have
\[
\norm{S(u)}^2
\le 3\norm{Au_t}^2+3\nu^2\norm{A^2u}^2
   +3\norm{B(Au,u)}^2
\le C\left(\norm{u_t}_{\dot V_2}^2
   +\nu^2\norm{u}_{\dot V_4}^2+M_3^4\right).
\]
Integrating over $\mathcal J$ and applying \eqref{eq:unit-bound},
we obtain
\[
\int_{\mathcal J}\norm{S(u(s))}^2\,ds
\le C\left((1+\nu^2)M_4^2(\ell+1)+M_3^4\ell\right).
\]
For $\ell\ge\tau_0>0$, we have
$\ell+1\le(1+\tau_0^{-1})\ell$. Hence
\[
\int_{\mathcal J}\norm{S(u(s))}^2\,ds
\le C\left((1+\nu^2)M_4^2(1+\tau_0^{-1})+M_3^4\right)\ell
= M_S^2\ell,
\]
where the constant in \eqref{eq:reference-constants} is chosen
sufficiently large. This proves \eqref{eq:S-control}.

For $t\ge r\ge t_1$, set $I_j=[t-j-1,t-j]\cap[r,t]$ for $j\ge0$.
Each nonempty $I_j$ is contained in a unit interval starting at or after $t_1$.
Applying Assumption~\ref{ass:reg} on these intervals, we obtain
\[
\begin{aligned}
\int_r^t e^{-\kappa(t-s)}\norm{S(u(s))}^2\,ds
&\le C\sum_{j=0}^{\infty}e^{-\kappa j}
\int_{I_j}\bigl(\norm{Au_t}^2+\nu^2\norm{A^2u}^2\bigr)\,ds
+\frac{CM_3^4}{\kappa}\\
&\le \frac{C(1+\nu^2)M_4^2}{1-e^{-\kappa}}
+\frac{CM_3^4}{\kappa}.
\end{aligned}
\]
\end{proof}

We have the following well-posedness result.

\begin{proposition}\label{prop:assim-wellposed}
Let $\beta_n\in[\alpha_0,\alpha_1]$, $\eta>0$, $w_n(t_n)\in V$,
$f\in L^\infty(0,\infty;H)$, and let $u$ be the reference solution. With $\beta_n$ fixed,
\eqref{eq:assim-system} has a unique global solution satisfying, for every
$T>t_n$,
\[
w_n\in C([t_n,T];V)\cap L^2(t_n,T;D(A)),\qquad
(w_n)_t\in L^2(t_n,T;H).
\]
\end{proposition}

\begin{proof}
Fix $n$ and $T>t_n$, and write $\beta=\beta_n$.
We use the Galerkin construction in \cite[Theorem~2]{MLAlphaParameter2025}. For $m\ge1$, let $\Pi_m$ be the orthogonal projections onto the spaces in
an increasing sequence of finite-dimensional Stokes eigenspaces
such that $P_NH\subset\Pi_mH$ and $\Pi_m\to I$ strongly in $H$ and $V$.
Let $w_m$ denote the corresponding Galerkin approximation, with
$w_m(t_n)=\Pi_mw_n(t_n)$.
Since $P_N$ is an orthogonal projection commuting with $A$,
\[
(T_\beta P_Nw_m,w_m)=E_\beta(P_Nw_m)\ge0.
\]
Testing the Galerkin equation with $w_m$, using
\eqref{eq:skew-symmetry} and Young's inequality, gives
\[
\frac{d}{dt}E_\beta(w_m)
+\nu\bigl(\norm{\nabla w_m}^2+\beta^2\norm{Aw_m}^2\bigr)
+\eta E_\beta(P_Nw_m)
\le \frac{F^2}{\nu\lambda_1}+\eta E_\beta(P_Nu).
\]
Thus, for every $\eta>0$, the approximations exist on $[t_n,T]$
and are uniformly bounded in
$L^\infty(t_n,T;V)\cap L^2(t_n,T;D(A))$.
Applying $T_\beta^{-1}$ to their equations and using
\eqref{eq:filter-inverse} and \eqref{eq:negative-bilinear}, we obtain
\[
\norm{(w_m)_t}
\le \nu\norm{Aw_m}
+C\alpha_0^{-2}\norm{T_\beta w_m}\norm{\nabla w_m}
+F+\eta\norm{P_N(w_m-u)}.
\]
Thus $(w_m)_t$ is bounded in $L^2(t_n,T;H)$.
By the Aubin--Lions theorem, after extracting a subsequence,
\[
w_m\to w_n\quad\text{in }L^2(t_n,T;V),\qquad
T_\beta w_m\rightharpoonup T_\beta w_n
\quad\text{in }L^2(t_n,T;H).
\]
The limit solves \eqref{eq:assim-system} with the given initial data. The embedding
$L^2(t_n,T;D(A))\cap H^1(t_n,T;H)
\hookrightarrow C([t_n,T];V)$ yields the continuity.

For uniqueness, let $e=w^{(1)}-w^{(2)}$ for two solutions with
identical data. Applying the interpolation estimates and Young's inequality, we have
\[
\big|\langle B(T_\beta e,w^{(1)}),e\rangle\big|
\le \frac{\nu\beta^2}{2}\norm{Ae}^2
+C\bigl(1+\norm{\nabla w^{(1)}}\norm{Aw^{(1)}}\bigr)E_\beta(e).
\]
Testing the difference equation with $e$  yields
\[
\frac{d}{dt}E_\beta(e)
\le C\bigl(1+\norm{\nabla w^{(1)}}\norm{Aw^{(1)}}\bigr)E_\beta(e).
\]
 Gronwall's
inequality implies $e=0$. Since $T>t_n$ is arbitrary,
the solution is global.
\end{proof}

For a sufficiently large domain-dependent constant $C$, we define
\begin{equation}\label{eq:C0}
C_0=C(1+\nu^{-3})\left(M_{\rm ref}^4\alpha_0^{-4}
+\alpha_1^2(1+M_3^4)\right).
\end{equation}
We require the observation resolution and nudging strength to satisfy
\begin{equation}\label{eq:N-lower}
\frac{2C_0}{\lambda_1N^2}\le\frac{\nu}{8},\qquad \eta>C_0.
\end{equation}
Using the quantities in \eqref{eq:reference-constants}, we also fix
\begin{equation}\label{eq:synchronization-constants}
\kappa=\frac{\nu\lambda_1}{8},\qquad
C_S^*=\frac{8}{\nu\lambda_1},\qquad
K_0=C_S^*K_S(\kappa),\qquad
K_A=\frac{4C_S^*M_S^2}{\nu\alpha_0^2}.
\end{equation}
The following estimate refines \cite{MLAlphaParameter2025} by retaining the higher-order dissipation and keeping the constants uniform over all assimilation intervals. This improvement is useful in the parameter update.

\begin{proposition}\label{prop:fixed}
Let $u$ be a solution of \eqref{eq:MLalpha-v} satisfying
Assumption~\ref{ass:reg} and $E_\alpha(u)\le M_{\rm ref}^2$
on $[t_n,t_{n+1}]$.
Let $w$ be a solution of \eqref{eq:assim-system} with fixed $\beta\in[\alpha_0,\alpha_1]$,
and set $g=w-u$. Under \eqref{eq:N-lower}, with constants
\eqref{eq:C0} and \eqref{eq:synchronization-constants}, we have
\begin{equation}\label{eq:fixed-differential}
\frac{d}{dt}E_\beta(g)+\kappa E_\beta(g)+\frac{\nu}{4}\norm{\nabla g}^2+\frac{\nu\beta^2}{4}\norm{Ag}^2
\le C_S^*\norm{S(u(t))}^2|\beta^2-\alpha^2|^2
\end{equation}
for a.e. $t\in[t_n,t_{n+1}]$. Moreover, for $t\in[t_n,t_{n+1}]$ and
$\mathcal J=[r,r+\ell]\subset[t_n,t_{n+1}]$ with $\ell\ge\tau_0>0$,
\begin{equation}\label{eq:fixed-estimates}
\begin{aligned}
E_\beta(g(t))&\le e^{-\kappa(t-t_n)}E_\beta(g(t_n))
+K_0|\beta^2-\alpha^2|^2,\\
\int_{\mathcal J} \norm{Ag(s)}^2ds&\le \frac{4}{\nu\beta^2}E_\beta(g(r))
+K_A|\beta^2-\alpha^2|^2\ell.
\end{aligned}
\end{equation}
The constants $\kappa$, $C_S^*$, $K_0$, and $K_A$ are independent of $N$ and $\eta$ once \eqref{eq:N-lower} holds.
\end{proposition}

\begin{proof}
Subtracting the reference equation from \eqref{eq:assim-functional}, we obtain
\begin{equation*}
T_\beta g_t+\nu AT_\beta g
+B(T_\beta g,u)+B(T_\beta w,g)
+\eta T_\beta P_Ng
=(\alpha^2-\beta^2)S(u).
\end{equation*}
Here we have used the identity $B(T_\beta w,w)-B(T_\beta u,u)=B(T_\beta g,u)+B(T_\beta w,g)$.
Pairing with $g$ and using \eqref{eq:skew-symmetry}, we have
\begin{align*}
\frac12\frac{d}{dt}E_\beta(g)&+\nu(\norm{\nabla g}^2+\beta^2\norm{Ag}^2)
+\eta\bigl(\norm{P_Ng}^2+\beta^2\norm{\nabla P_Ng}^2\bigr)\\
&\le |\langle B(g,u),g\rangle|+\beta^2|\langle B(Ag,u),g\rangle|
+|\beta^2-\alpha^2|\,|(S(u),g)|.
\end{align*}
Using the interpolation inequalities, Young's inequality, and $\norm{\nabla u}\le M_{\rm ref}/\alpha_0$,
$\norm{\nabla u}_{L^\infty}\le CM_3$, and $\beta\le\alpha_1$, we have
\[
\begin{aligned}
|\langle B(g,u),g\rangle|&\le \frac{\nu}{16}\norm{\nabla g}^2
+C\nu^{-3}\norm{\nabla u}^4\norm{g}^2,\\
\beta^2|\langle B(Ag,u),g\rangle|&\le \frac{\nu\beta^2}{16}\norm{Ag}^2
+C\nu^{-1}\beta^2 M_3^2\norm{g}^2.
\end{aligned}
\]
The Poincar\'e and Young inequalities give
\begin{equation*}
|\beta^2-\alpha^2|\,|(S(u),g)|\le \lambda_1^{-1/2}|\beta^2-\alpha^2|\,\norm{S(u(t))}\norm{\nabla g}
\le \frac{\nu}{16}\norm{\nabla g}^2+\frac{4}{\nu\lambda_1}\norm{S(u(t))}^2|\beta^2-\alpha^2|^2.
\end{equation*}
Choosing $C_0$ large enough, we then obtain
\begin{align*}
\frac{d}{dt}E_\beta(g)&+\frac{\nu}{2}\norm{\nabla g}^2+\frac{\nu\beta^2}{2}\norm{Ag}^2+2\eta(\norm{P_Ng}^2+\beta^2\norm{\nabla P_Ng}^2)\\
&\le C_0\norm{g}^2+C_S^*\norm{S(u(t))}^2|\beta^2-\alpha^2|^2.
\end{align*}
Using the orthogonality of the low and high Fourier modes together with \eqref{eq:projection-properties}, we have
\[
C_0\norm{g}^2\le 2C_0\norm{P_Ng}^2+\frac{2C_0}{\lambda_1N^2}\norm{\nabla g}^2.
\]
By \eqref{eq:N-lower}, the low-mode contribution is absorbed by the feedback term, and the high-mode contribution is absorbed by part of the viscous dissipation. Consequently,
\[
\frac{d}{dt}E_\beta(g)+\frac{3\nu}{8}\norm{\nabla g}^2+\frac{\nu\beta^2}{2}\norm{Ag}^2
\le C_S^*\norm{S(u(t))}^2|\beta^2-\alpha^2|^2.
\]
From \eqref{eq:poincare}, we have
\[
\frac{\nu}{8}\norm{\nabla g}^2+\frac{\nu\beta^2}{4}\norm{Ag}^2
\ge\frac{\nu\lambda_1}{8}
\bigl(\norm{g}^2+\beta^2\norm{\nabla g}^2\bigr).
\]
This proves \eqref{eq:fixed-differential}.
Multiplying by the exponential integrating factor and applying \eqref{eq:S-weighted}, we obtain the first bound in
\eqref{eq:fixed-estimates}. Integrating the differential inequality over
$\mathcal J$ and applying \eqref{eq:S-control}, we obtain
\[
\frac{\nu\beta^2}{4}\int_{\mathcal J}\norm{Ag(s)}^2ds\le E_\beta(g(r))+C_S^*M_S^2|\beta^2-\alpha^2|^2\ell.
\]
Since $\beta\ge\alpha_0$, we obtain the second bound in \eqref{eq:fixed-estimates}.
\end{proof}

\begin{remark}
For the synchronization estimates above, once $N$ satisfies
\eqref{eq:N-lower}, any fixed $\eta>C_0$ is allowed.
Unlike condition~(20) in \cite{BardinaParameter2026},
these estimates impose no upper bound on $\eta$.
\end{remark}

\section{Synchronization estimates}\label{sec:convergence}

We use the synchronization and dissipation estimates to control the error introduced by replacing the unobserved part of the reference solution.
The following identity relates this error to the parameter update.

\begin{lemma}\label{lem:update-id}
Under Assumption~\ref{ass:reg}, let $w_n$ be given by
Proposition~\ref{prop:assim-wellposed}. For a.e. $t\in I_n$, there is a residual
$\mathcal R_N^{(n)}$ satisfying
\begin{equation}\label{eq:update-id}
\begin{aligned}
Y_N^{(n)}(t)&=(\alpha^2-\beta_n^2)D_N^{(n)}(t)+\mathcal R_N^{(n)}(t),\\
\norm{\mathcal R_N^{(n)}(t)}
&\le C_R\bigl(1+E_{\beta_n}(g_n(t))^{1/2}\bigr)
\bigl(E_{\beta_n}(g_n(t))^{1/2}+\norm{AQ_Ng_n(t)}\bigr),
\end{aligned}
\end{equation}
where $C_R$ depends on $\nu,\alpha_0,\alpha_1,M_3,N$, but is independent of $n$.
\end{lemma}

\begin{proof}
Applying the low-mode projection to \eqref{eq:key-identity} and using the definitions in \eqref{eq:Dn} and \eqref{eq:Yn}, we obtain the identity in \eqref{eq:update-id}, with residual
\[
\mathcal R_N^{(n)}
=P_N\bigl\{B(T_\alpha u,u)-B(T_\alpha u_N^{*,n},u_N^{*,n})\bigr\}.
\]
To estimate this residual, we write $h=Q_Ng_n$, so that $u_N^{*,n}=u+h$.

For $\phi\in P_NH$ with $\norm{\phi}=1$, we have
$\norm{\nabla\phi}_{L^\infty}\le c_N$, where $c_N$ depends
only on $N$ and $\Omega$.
By bilinearity and integration by parts, using
$\nabla\cdot u=\nabla\cdot h=0$, we obtain
\[
\begin{aligned}
\big|\langle B(u,u)-B(u+h,u+h),\phi\rangle\big|
&\le \bigl[c_N(\norm{u}+\norm{h})
+\norm{\nabla u}_{L^\infty}\bigr]\norm{h},\\
\big|\langle B(Au,u)-B(A(u+h),u+h),\phi\rangle\big|
&\le c_N(\norm{Au}+\norm{Ah})\norm{h}
+\norm{\nabla u}_{L^\infty}\norm{Ah}.
\end{aligned}
\]
Combining these estimates and taking the supremum over unit test fields $\phi\in P_NH$, we obtain
\[
\norm{\mathcal R_N^{(n)}}\le C(N,\alpha_1,\Omega)
\bigl(\norm{u}+\norm{Au}+\norm{\nabla u}_{L^\infty}+\norm{h}\bigr)(\norm{h}+\norm{Ah}).
\]
Finally, we use \eqref{eq:reg1}, the embedding $\dot V_3\hookrightarrow W^{1,\infty}$,
the inequality $\norm{h}\le E_{\beta_n}(g_n)^{1/2}$, and the identity $\norm{Ah}=\norm{AQ_Ng_n}$.
These bounds yield \eqref{eq:update-id}, with a constant uniform in $n,L$.
\end{proof}

For $\delta_n>0$, substituting \eqref{eq:update-id} into \eqref{eq:update} gives
\[
\gamma_{n+1}-\alpha^2=
\delta_n^{-1}\int_{I_n}(\mathcal R_N^{(n)},D_N^{(n)})\,ds.
\]
Applying the Cauchy--Schwarz inequality, we obtain
\begin{equation}\label{eq:one-step-update}
\begin{aligned}
|\beta_{n+1}^2-\alpha^2|
&\le\frac{1}{\delta_n}\int_{I_n}
\norm{\mathcal R_N^{(n)}}\,\norm{D_N^{(n)}}\,ds\\
&\le\frac{\norm{\mathcal R_N^{(n)}}_{L^2(I_n;H)}}{\delta_n^{1/2}}
=\frac{\tau_n^{-1/2}\norm{\mathcal R_N^{(n)}}_{L^2(I_n;H)}}
{(\delta_n/\tau_n)^{1/2}}.
\end{aligned}
\end{equation}
With $C_R$ as in Lemma~\ref{lem:update-id}, we fix the following constants before choosing the relaxation length $L$:
\[
P_*=\alpha_1^2-\alpha_0^2,\qquad E_*=1+K_0P_*^2,
\qquad C_A=\frac{2}{\alpha_0\sqrt{\nu\tau_0}}.
\]
We then define
\begin{equation}\label{eq:Crec}
C_{\rm rec}=C_R(1+\sqrt{E_*})
\max\left\{1+C_A,(1+C_A)\sqrt{K_0}+\sqrt{K_A}\right\}.
\end{equation}
The constants $K_0,K_A,E_*,C_{\rm rec}$ are independent of $L$ and $w_0$.
The remaining assumptions require successful update windows
and a sufficiently strong signal relative to the bound on
$\mathcal R_N^{(n)}$.

\begin{assumption}\label{ass:conv}
\begin{enumerate}[label=\textup{(H\arabic*)},leftmargin=3.0em]
\item $\alpha\in[\alpha_0,\alpha_1]$, with
$\alpha_0>0$. The reference solution satisfies Assumption~\ref{ass:reg}
and the initial energy bound used in \eqref{eq:Mref}.
\item  $N$ and $\eta$ remain fixed and satisfy
\eqref{eq:N-lower}, with $C_0$ given by \eqref{eq:C0}.
\item Fix the trial lengths $0<\ell_1<\cdots<\ell_J\le\tau_{\max}<\infty$
and threshold $\zeta>0$ independently of $L$.
There exists $L_{\rm sig}>0$ such that, for every $L\ge L_{\rm sig}$,
Algorithm~\ref{alg:recovery}, run without a terminal time,
accepts a trial at every iteration. Hence
$\tau_n\in\{\ell_1,\ldots,\ell_J\}$ and
\begin{equation}\label{eq:zeta}
\delta_n/\tau_n\ge \zeta>0\qquad\text{for every }n.
\end{equation}
\item  With $C_{\rm rec}$ from
\eqref{eq:Crec},
\begin{equation}\label{eq:contraction-coeff}
C_{\rm rec}\zeta^{-1/2}<1.
\end{equation}
\end{enumerate}
\end{assumption}

\begin{remark}\label{rem:assumption-feasibility}
Conditions \textup{(H3)}--\textup{(H4)} are additional assumptions:
\textup{(H3)} requires a successful update window at every iteration,
while \textup{(H4)} requires $\zeta>C_{\rm rec}^2$.
Neither follows from Assumption~\ref{ass:reg}.
Unlike \cite{BardinaParameter2026}, our bound on
$\mathcal R_N^{(n)}$ has no explicit decay in $N$; the constant $C_R$
contains $c_N$, which may grow with $N$.
Thus increasing $N$ does not necessarily improve \textup{(H4)}.
\end{remark}

\begin{theorem}\label{thm:main}
Under Assumption~\ref{ass:conv}, generate $w_n,\beta_n$ by
\eqref{eq:assim-system}, \eqref{eq:update} and \eqref{eq:clipping}, from
$w_0\in V$, $\beta_1\in[\alpha_0,\alpha_1]$, with $\widehat t_n=t_n+L$.
There exists $L_0\ge L_{\rm sig}$, which may depend on $w_0$, such that every fixed $L\ge L_0$ admits
$C_L>0$ and $\theta_L\in(0,1)$ satisfying
\begin{equation*}
|\beta_n^2-\alpha^2|\le C_L\theta_L^{n-1},\qquad
E_{\beta_n}\big(w_n(t_n)-u(t_n)\big)^{1/2}\le C_L\theta_L^{n-1}.
\end{equation*}
Moreover, there exist constants $\widetilde C_L,\sigma_L>0$ such that
\[
|\beta_n^2-\alpha^2|
+E_{\beta_n}\bigl(w_n(t_n)-u(t_n)\bigr)^{1/2}
\le \widetilde C_L e^{-\sigma_L(t_n-t_1)},
\qquad n\ge1,
\]
where $\sigma_L=|\log\theta_L|/(L+\tau_{\max})$.
In particular, $\beta_n\to\alpha$ and $w_n(t_n)\to u(t_n)$ in $V$ as $n\to\infty$. 
\end{theorem}

\begin{proof}
Set
$a_n:=|\beta_n^2-\alpha^2|$,
$b_n:=E_{\beta_n}(w_n(t_n)-u(t_n))^{1/2}$, and 
\[
B_*:=\left(\norm{w_0}^2+\alpha_1^2\norm{\nabla w_0}^2\right)^{1/2}+\frac{\alpha_1}{\alpha_0}M_{\rm ref}.
\]
Since $\alpha,\beta_1\in[\alpha_0,\alpha_1]$, the triangle
inequality and \eqref{eq:Mref} give
\[
\begin{aligned}
b_1
&\le E_{\beta_1}(w_0)^{1/2}
   +E_{\beta_1}(u(t_1))^{1/2}\\
&\le E_{\alpha_1}(w_0)^{1/2}
   +\frac{\alpha_1}{\alpha_0}E_\alpha(u(t_1))^{1/2}
\le B_*.
\end{aligned}
\]
Moreover, (\ref{eq:clipping}) ensures
$\beta_n\in[\alpha_0,\alpha_1]$, so
\[
a_n=|\beta_n^2-\alpha^2|
\le\alpha_1^2-\alpha_0^2=P_*.
\]
Applying the first estimate in \eqref{eq:fixed-estimates} at the end of the assimilation interval, we obtain
\begin{equation*}
E_{\beta_n}(g_n(t_{n+1}))^{1/2}\le e^{-\kappa(t_{n+1}-t_n)/2}b_n+K_0^{1/2}a_n.
\end{equation*}
We set $C_{\rm eq}=\alpha_1/\alpha_0$ and $C_1=C_{\rm eq}K_0^{1/2}$. Using $E_{\beta_{n+1}}\le C_{\rm eq}^2E_{\beta_n}$ and $t_{n+1}-t_n\ge L$, we conclude that
\begin{equation*}
b_{n+1}\le q_Lb_n+C_1a_n,\qquad q_L:=C_{\rm eq}e^{-\kappa L/2}.
\end{equation*}
Choose $L_{\rm pre}\ge L_{\rm sig}$ sufficiently large that $q_L\le1/2$ for every $L\ge L_{\rm pre}$.
Iterating the recurrence and using $a_n\le P_*$ gives
\[
b_n\le q_L^{n-1}B_*+\frac{C_1P_*(1-q_L^{n-1})}{1-q_L}
\le B_{\rm bd}:=\max\{B_*,2C_1P_*\}.
\]
We increase $L_{\rm pre}$, if necessary, so that $e^{-\kappa L}B_{\rm bd}^2\le1$.
The first estimate in \eqref{eq:fixed-estimates} then gives, on every update interval $I_n$,
\[
E_{\beta_n}(g_n(s))\le e^{-\kappa L}B_{\rm bd}^2+K_0P_*^2\le E_*.
\]
This bound makes the coefficient in the residual estimate uniform over the update intervals. Denoting $E=E_{\beta_n}(g_n)\le E_*$ and using $\norm{AQ_Ng_n}\le\norm{Ag_n}$, we combine
\eqref{eq:one-step-update}, Lemma~\ref{lem:update-id}, and \eqref{eq:zeta}, and apply
Minkowski's inequality to obtain
\begin{equation*}
a_{n+1}\le C_R(1+\sqrt{E_*})\zeta^{-1/2}\left[\left(\frac1{\tau_n}\int_{I_n}E(s)ds\right)^{1/2}+\left(\frac1{\tau_n}\int_{I_n}\norm{Ag_n(s)}^2ds\right)^{1/2}\right].
\end{equation*}
Applying \eqref{eq:fixed-estimates} over the update interval, since $\tau_n\ge\tau_0$ and
$E(\widehat t_n)^{1/2}\le e^{-\kappa L/2}b_n+\sqrt{K_0}a_n$, we have
\[
\begin{aligned}
\left(\frac1{\tau_n}\int_{I_n}E(s)\,ds\right)^{1/2}
&\le e^{-\kappa L/2}b_n+\sqrt{K_0}a_n,\\
\left(\frac1{\tau_n}\int_{I_n}\norm{Ag_n(s)}^2\,ds\right)^{1/2}
&\le C_Ae^{-\kappa L/2}b_n+(C_A\sqrt{K_0}+\sqrt{K_A})a_n.
\end{aligned}
\]
Substituting these bounds into the preceding estimate for $a_{n+1}$
and using \eqref{eq:Crec}, we obtain
\begin{equation*}
a_{n+1}\le C_{\rm rec}\zeta^{-1/2}a_n+C_{\rm rec}\zeta^{-1/2}e^{-\kappa L/2}b_n.
\end{equation*}
Combining this estimate with
$b_{n+1}\le C_{\rm eq}e^{-\kappa L/2}b_n+C_1a_n$,
we obtain the recurrence 
\[
\mathcal{T}_L:=\begin{pmatrix}
C_{\rm rec}\zeta^{-1/2} & C_{\rm rec}\zeta^{-1/2}e^{-\kappa L/2}\\
C_1 & C_{\rm eq}e^{-\kappa L/2}
\end{pmatrix},\qquad
\binom{a_{n+1}}{b_{n+1}}\le \mathcal{T}_L\binom{a_n}{b_n}.
\]
Let $\rho(\mathcal T_L)$ denote the spectral radius of $\mathcal T_L$.
As $L\to\infty$,
\[
\mathcal T_L\longrightarrow
\begin{pmatrix}
C_{\rm rec}\zeta^{-1/2} & 0\\
C_1 & 0
\end{pmatrix}
.
\]
The limiting matrix has eigenvalues
$C_{\rm rec}\zeta^{-1/2}$ and $0$. Hence, by continuity
of the spectral radius and \textup{(H4)},
\[
\rho(\mathcal T_L)\longrightarrow
C_{\rm rec}\zeta^{-1/2}<1.
\]
Thus $\rho(\mathcal T_L)<1$ for all sufficiently large $L$.
Choose $L_0\ge L_{\rm pre}$ so that this holds for every $L\ge L_0$.
By \textup{(H3)} and the nonnegativity of $\mathcal T_L$, induction yields
\[
\binom{a_n}{b_n}\le\mathcal T_L^{n-1}\binom{a_1}{b_1}.
\]
For $\rho(\mathcal T_L)<\theta_L<1$, the standard matrix-power estimate gives, for some $C_{\theta,L}>0$,
\[
0\leq a_n+b_n\le C_{\theta,L}\theta_L^{n-1}(a_1+b_1).
\]
This proves the geometric bounds. Moreover, since $0<\theta_L<1$, the preceding estimate gives
$a_n\to0$ and $b_n\to0$. By their definitions and
$\alpha,\beta_n\ge\alpha_0>0$, we obtain
\[
|\beta_n-\alpha|\le\frac{a_n}{2\alpha_0}\to0,\qquad
\norm{\nabla(w_n(t_n)-u(t_n))}\le\frac{b_n}{\alpha_0}\to0.
\]
Finally, since
$n-1\ge(t_n-t_1)/(L+\tau_{\max})$, we have
\[
\theta_L^{n-1}
=e^{(n-1)\log\theta_L}
\le \exp\left(-\frac{|\log\theta_L|}{L+\tau_{\max}}
(t_n-t_1)\right).
\]
This gives exponential decay in $t_n-t_1$ at the update times.
For a fixed geometric factor, increasing $L$ decreases this rate,
so a longer relaxation interval does not by itself guarantee faster convergence.
\end{proof}

\section{Conclusions}\label{sec:conclusion}

In this study, we have constructed a recursive algorithm for recovering the length-scale parameter in the modified Leray-$\alpha$ model. The synchronization estimate controls the velocity error,
while the integrated dissipation estimate supplies the higher-order control needed for the analysis.
With fixed $N,\eta$, a sufficiently large fixed relaxation length,
uniform reference regularity, successful update windows, and the condition $C_{\rm rec}\zeta^{-1/2}<1$,
the parameter and velocity errors converge geometrically at the update times. If the update windows have a uniform upper bound on their lengths,
the errors also decay exponentially in time at the update times.

The convergence result remains conditional on the signal hypotheses. Deriving these conditions from the reference trajectory remains open.
Martinez, Strikwerda, and Wan \cite{MartinezStrikwerdaWan2026} study excitation
and viscosity identifiability for the two-dimensional Navier--Stokes equations
using spectral observations and a nudging-based sensitivity analysis.
A natural direction for future work is to investigate whether a similar analysis
can provide verifiable sufficient conditions for \textup{(H3)}--\textup{(H4)}. Moreover, the present analysis assumes exact observations available continuously
in time. Extensions to discrete-time observations and noisy data
are left for future work.

\section*{Acknowledgments}

Jing Tian’s work was partially supported by the NSF LEAPS-MPS Grant \#2316894.

\end{document}